\documentclass[11.5pt]{paper}

\usepackage{hyperref}
\usepackage{graphicx}
\usepackage{amsmath}
\usepackage{xcolor}
\usepackage{amssymb}
\usepackage{amsfonts}
\usepackage{marginnote}
\usepackage{stmaryrd}
\usepackage{harpoon}

\def\ddd{n}

\def\omegen{{\widetilde\Omega}}

\def\vt{\vartheta}
\def\vtn{\vartheta_M}
\def\vtnb{{\vartheta}_M^*}
\def\vta{\vartheta^*}
\def\ua{u^\ast}
\def\unt{\widetilde{u}_N}
\def\vtnt{\widetilde{\vartheta}_M}

\newcommand{\mcC}{\mathcal{C}}
\newcommand{\mcD}{\mathcal{D}}
\newcommand{\mcG}{{\mathcal{D}^\ast}}

\newcommand{\mcV}{\mathcal{V}}

\newcommand{\mcL}{\mathcal{L}}

\newcommand{\mbR}{\mathbb{R}}

\newcommand{\mbRn}{{\mathbb{R}^n}}

\def\omgsc{{\omgs\cup\omgc}}
\newcommand{\omg}{{\Omega}}
\newcommand{\omgs}{{\Omega}}
\newcommand{\omgc}{{\Omega_{\mathcal I}}}

\def \alphab{{\boldsymbol\alpha}}

\def \nub{{\boldsymbol \nu}}

\def \bthe{{\boldsymbol \Theta}}

\def \xb{\mathbf{x}}

\def \yb{\mathbf{y}}

\newtheorem{thm}{Theorem}[section]
\newtheorem{lemma}{Lemma}[section]

\begin{document}

\thispagestyle{empty}
\begin{center}
\textbf{\Large Identification of the diffusion parameter in\\ 
\smallskip
nonlocal steady diffusion problems}\footnote{This work was supported in part by the US National Science Foundation grant DMS-1315259.}

\bigskip M. D'Elia\footnote{Currently at Sandia National Laboratories, NM. Email: {\tt mdelia@sandia.gov}. Sandia National Laboratories is a multi program laboratory managed and operated by Sandia Corporation, a wholly owned subsidiary of Lockheed Martin Corporation, for the U.S. Department of Energy's National Nuclear Security Administration under contract DE-AC04-94AL85000.} and M. Gunzburger 

\vspace{.2cm} Florida State University, 400 Dirac Science Library, Tallahassee, FL, 32306, USA
\end{center}

\begin{abstract}
The problem of identifying the diffusion parameter appearing in a nonlocal steady diffusion equation is considered. The identification problem is formulated as an optimal control problem having a matching functional as the objective of the control and the parameter function as the control variable. The analysis makes use of a nonlocal vector calculus that allows one to define a variational formulation of the nonlocal problem. In a manner analogous to the local partial differential equations counterpart, we demonstrate, for certain kernel functions, the existence of at least one optimal solution in the space of admissible parameters. We introduce a Galerkin finite element discretization of the optimal control problem and derive a priori error estimates for the approximate state and control variables. Using one-dimensional numerical experiments, we illustrate the theoretical results and show that by using nonlocal models it is possible to estimate non-smooth and discontinuous diffusion parameters.
\end{abstract}

\vspace{.2cm}

\section{Introduction}
Nonlocal models are currently used in a large variety of applications, including continuum mechanics \cite{chgu:11,sill:00}, graph theory \cite{lova:06}, nonlocal wave equations \cite{weab:05}, and jump processes \cite{bbck:09,bakk:10,burc:11}. Our interest is in nonlocal diffusion operators that arise in many fields where the phenomena cannot be modeled accurately by the standard classical diffusion equation. Among those applications we have image analyses \cite{bucm:10,gilboa:595,gilboa:1005,lzob:10}, machine learning \cite{robd:10}, nonlocal Dirichlet forms \cite{appl:04}, kinetic equations \cite{bass:84,limi:10}, phase transitions \cite{bach:99,fife:03}, and nonlocal heat conduction \cite{bodu:09}. 

The principal difference between nonlocal models and the classical (local) partial differential equations (PDEs) is in how two or more domains interact. In the local case interactions only occur due to contact, whereas in the nonlocal case they can occur at distance. Consider an open bounded domain $\omgs$ in $\mbRn$. For $u(\xb)\colon \Omega \to \mbR$, the action of the nonlocal diffusion operator $\mcL$ on the function $u(\xb)$ is defined as
\begin{equation}\label{nldo}
   \mcL u(\xb) := 2\int_{\mbR^n} \big(u(\yb)-u(\xb)\big) \, \vt(\xb,\yb)\,\gamma (\xb, \yb )\,d\yb \qquad  \forall\,\xb \in \omgs \subseteq \mbRn,
\end{equation}
where the kernel $\gamma (\xb, \yb )\colon\Omega\times\Omega\to\mbR$ is a non-negative symmetric mapping (i.e., $\gamma (\xb, \yb )=\gamma (\yb, \xb )$), and the diffusion parameter $\vt(\xb, \yb )\colon\Omega\times\Omega\to\mbR$ is a positive function. We are interested in the nonlocal steady-state diffusion equation 
\begin{equation}\label{nlde}
\left\{\begin{array}{ll}
-\mcL u = f &\qquad \mbox{on $\omgs$} \\
u  = g &\qquad \mbox{on $\omgc$},
\end{array}\right.
\end{equation} 
where the equality constraint (the nonlocal counterpart of a Dirichlet boundary condition for PDEs) acts on an interaction volume $\omgc$ that is disjoint from $\omgs$. Nonlocal diffusion problems such as \eqref{nlde} have been analyzed in the recent works \cite{akme:10,dglz:11,Du10,gule:10}, and techniques for an accurate numerical solution have been developed and applied to diverse applications \cite{chgu:11,Aksoylu,bule:11a,bule:12,Du11,spgl09,zhdu:10}. However, these mathematical models are not exact; parameters such as volume constraint data, diffusivity coefficients, and source terms are often unknown or subject to uncertainty. This fact affects the quality of the computational results and their reliability, making necessary the development of techniques for the identification of model parameters. In the classical local case, i.e., for PDEs, a widely used approach is to formulate the problem as an optimal control problem having the unknown parameters acting as control variables; here, we follow the same approach for nonlocal problems.  Control problems for nonlocal diffusion equations have already been addressed in \cite{degu:13} where, for the same nonlocal model \eqref{nlde}, the control variable is the source term $f$; this results in a {\it linear} optimization problem constrained by the nonlocal equation. The case treated in this work is more complex (as its local counterpart is) because of its {\it nonlinearity} and thus requires more sophisticated techniques for analyzing the well-posedness and designing finite dimensional approximations. The estimation problem must rely on some additional a priori information, e.g., a target or reference function for the function of interest $u$. In our context, we formulate the parameter identification as the problem of finding $u$ and $\vt$ such that $u$ is as close as possible to a target function $\widehat{u}$ under the constraint that $u$ and $\vt$ satisfy the nonlocal problem \eqref{nlde}. 

The main contribution of this work is to show that we can mimic the approaches of the classical control theory for PDEs and to show that by using nonlocal models one can obtain accurate estimates of non-smooth and discontinuous diffusion parameters, a case that often arises in practice. 

Our analysis is based on the nonlocal vector calculus introduced in \cite{Du10} which is exploited to define a variational formulation of the nonlocal equation and for studying the well-posedness of the control problem. The use of the nonlocal calculus is crucial, as is the classical calculus in the PDE setting, to the analyses. In fact, it helps one avoid more cumbersome direct approaches to deriving the results we present in this work.

 Relevant aspects of the nonlocal calculus are reviewed in Section \ref{nlvc}. In Section \ref{sec:control} we formulate the identification problem as a control problem and we demonstrate the existence of at least one optimal solution. In Section \ref{fin_dim_approx} we study finite dimensional approximations of the optimal state and control and, for finite element approximations, we derive a priori error estimates. In Section \ref{num_tests} we present numerical tests conducted on one-dimensional problems and illustrate the theoretical results of Section \ref{fin_dim_approx}. These numerical results serve as a base for two- and three-dimensional simulations and show that by using nonlocal models it is possible to accurately approximate non-smooth parameter functions. In Section \ref{conclusion} we provide a few concluding remarks.

\section{A nonlocal vector calculus and volume-constrained problems}\label{nlvc}
For the sake of completeness, in this section we present the basic notions of the nonlocal vector calculus, recall some theoretical results that are useful to us in this work, and introduce nonlocal diffusion problems.
\subsection{Notation}
A detailed introduction to the nonlocal vector calculus can be found in \cite{Du10}. Here, we limit the discussion to the tools that we use throughout the paper. For the vector mappings $\nub(\xb,\yb), \alphab(\xb,\yb) \colon \mbRn\times\mbRn\to \mbRn$, with $\alphab$ antisymmetric (i.e., $\alpha (\xb, \yb )=-\alpha (\yb, \xb )$), the action of the nonlocal divergence operator $\mcD\colon \mbRn \to \mbR$ on $\nub$ is defined as
\begin{subequations}\label{ndivgrad}
\begin{equation}\label{ndiv}
 \mcD\big(\nub\big)(\xb) := \int_{\mbRn} \big(\nub(\xb,\yb)+\nub(\yb,\xb)\big)\cdot\alphab(\xb,\yb)\,d\yb\qquad
        \mbox{for $\xb\in\mbRn$}.
\end{equation}

The action of the adjoint operator $\mcG\colon  \mbRn\times\mbRn\to\mbRn$ on a mapping $u(\xb)\colon \mbRn$ $\to\mbR$ is given by
\begin{equation}\label{ngra}
\mcG\big(u\big)(\xb,\yb) = -\big(u(\yb)-u(\xb)\big)  \alphab(\xb,\yb) \qquad\mbox{for $\xb,\yb\in\mbRn$}.
\end{equation}
\end{subequations}
Thus, $-\mcG$ defines a nonlocal gradient operator.
We define the nonlocal diffusion operator $\mcL\colon \mbRn \to \mbR$ as the composition of the nonlocal divergence and gradient operators, i.e.  $\mcL u  := -\mcD\big(\vt \,\mcG u)$, where the diffusion parameter $\vt(\xb,\yb)$ is a positive symmetric function that maps $\mbRn\times\mbRn$ into $\mbR$. Then,  
\begin{equation}
\hspace{-2cm}\mcL u (\xb) := -\mcD\big(\vt \,\mcG u)(\xb)= 2\int_{\mbRn}\big(u(\yb)-u(\xb)\big)\,\vt(\xb,\yb) \,\gamma(\xb,\yb) \,d\yb\quad \hbox{for}\;\;\xb\in\mbRn
\end{equation}
which is exactly the operator we introduced in \eqref{nldo}. Thus, the nonlocal calculus allows us to express the nonlocal diffusion operator \eqref{nldo} as a composition of a nonlocal divergence operator and a nonlocal gradient operator. Here the symmetric kernel $\gamma$ is such that $\gamma(\xb,\yb):=\alphab(\xb,\yb)\cdot \alphab(\xb,\yb)$\footnote{In \cite{dglz:11} the diffusivity is defined as a second order symmetric positive definite tensor $\bthe$ and the kernel is defined as $\gamma:=\alphab \cdot \left(\bthe\alphab\right)$. Here for simplifying the analysis of the identification problem we consider a scalar diffusivity and we do not include it in the definition of the kernel.}.

Given an open subset $\omgs\subset\mbRn$, we define the interaction domain corresponding to $\omgs$ as
$$
   \omgc := \{ \yb\in\mbRn\setminus\omgs \quad\mbox{such that}\quad \alphab(\xb,\yb)\ne{\bf 0}\quad \mbox{for $\xb\in\omgs$}\}.
$$
Thus, $\omgc$ consists of those points outside of $\omgs$ that interact with points in $\omgs$.
\subsection{The kernel}\label{sec:kernel}
We assume that the domains $\omgs$, $\omgc$, and $\omgsc$ are  bounded with piecewise smooth boundary and satisfy the interior cone condition. We also assume that the symmetric kernel satisfies
\begin{equation}\label{gamma-conds}
\left\{\begin{array}{ll}
   \gamma(\xb,\yb) \geq  0 \quad &\forall\, \yb\in B_\varepsilon(\xb)\\[2mm]
   \gamma(\xb,\yb) \ge   \gamma_0 >  0 \quad &\forall\, \yb\in B_{\varepsilon/2}(\xb)\\[2mm]
   \gamma(\xb,\yb) = 0   \quad &\forall\, \yb\in (\omgsc) \setminus B_\varepsilon(\xb)
\end{array}\right.
\end{equation}
for all $\xb\in\omgsc$, where $\gamma_0$ and $\varepsilon$ are given positive constants and $B_\varepsilon({\xb}) := \{ \yb \in\omgsc  \colon |\yb-\xb|\le \varepsilon \}$; thus, nonlocal interactions are limited to a ball of radius $\varepsilon$ which is referred to as the interaction radius. This implies that 
\begin{equation}\label{omgie}
  \omgc = \{ \yb\in \mbR^\ddd\setminus\omg \,\,\,\colon\,\,\, |\yb-\xb|<\varepsilon \mbox{ for $\xb\in\omg$}\}.
\end{equation}

In \cite{dglz:11} (and also in \cite{akme:10, akpa:11,amrt:10}) several choices for the kernel $\gamma$ are considered and analyzed. Here, for the sake of brevity, we limit ourselves to two kernel classes. The results presented in this paper can be generalized to the other kernel functions considered in the above cited papers.

\vspace{.2cm}\noindent{\bf Case 1}. We further assume that there exist $s\in (0,1)$ and positive constants $\gamma_1$ and $\gamma_2$ such that, for all $\xb\in\omgsc$,
\begin{equation}\label{case1}
\frac{\gamma_1}{|\yb-\xb|^{\ddd+2s}} \leq \gamma(\xb,\yb) \leq \frac{\gamma_2}{|\yb-\xb|^{\ddd+2s}} \qquad \mbox{for $\yb\in B_\varepsilon({\xb})$}.
\end{equation}
An example is given by
$$
   \gamma(\xb,\yb) = \frac{\sigma(\xb,\yb)}{|\yb-\xb|^{\ddd+2s}}
$$
with $\sigma(\xb,\yb)$ symmetric and bounded from above and below by positive constants. 

\vspace{.2cm}\noindent{\bf Case 2}. In addition to \eqref{gamma-conds}, we assume that there exist positive constants $\gamma_3$ and $\gamma_4$ such that, for all $\xb\in\omgsc$,
\begin{equation}\label{case3}
\frac{\gamma_3}{|\yb-\xb|^{\ddd}} \leq \gamma(\xb,\yb) \leq \frac{\gamma_4}{|\yb-\xb|^{\ddd}} \qquad \mbox{for $\yb\in B_\varepsilon({\xb})$}.
\end{equation}
An example for this case is given by
$$
   \gamma(\xb,\yb) = \frac{\xi(\xb,\yb)}{|\yb-\xb|^{\ddd}}
$$
with $\xi(\xb,\yb)$ symmetric and bounded from above and below by positive constants. 
\subsection{Equivalence of spaces}\label{sec:equivsp}
We define the nonlocal energy semi-norm, nonlocal energy space, and nonlocal volume-constrained energy space by
\begin{subequations}
\begin{equation}
|||v|||^2 := \int_\omgsc\int_{\omgsc}\mcG(v)(\xb,\yb )\cdot\mcG(v)(\xb,\yb )\,d\yb \, d\xb \label{energynorm}
\end{equation}
\begin{equation}
V(\omgsc)  := \left\{ v  \in L^2(\omgsc) \,\,:\,\, |||v||| < \infty \right\}\qquad\quad\; \label{vspace}
\end{equation}
\begin{equation}
V_c(\omgsc)  := \left\{v\in V(\omgsc) \,\,:\,\, v=0\;{\rm on}\;\omgc\right\}.\qquad\;\;\label{vcspace}
\end{equation}
\end{subequations}
In \cite{dglz:11}, it is shown that, for kernels satisfying \eqref{gamma-conds} and \eqref{case1}, the nonlocal energy space $V(\omgsc)$ is equivalent to the fractional-order Sobolev space $H^s(\omgsc)$\footnote{For $s\in(0,1)$ and for a general domain $\omegen\in\mbR^\ddd$, let
$$
|v|_{H^s(\omegen)}^2 :=
\int_\omegen\int_\omegen\frac{\big(v(\yb)-v(\xb)\big)^2}{|\yb-\xb|^{\ddd+2s}}\,d\yb d\xb.
$$
Then, the space $H^s(\omegen)$ is defined by \cite{Adams}
$
H^s(\omegen) := \left\{v\in L^2(\omegen) : \|v\|_{L^2(\omegen)} +
|v|_{H^s(\omegen)}<\infty\right\}.
$
}.
This implies that $V_c(\omgsc)$ is a Hilbert space equipped with the norm $|||\cdot|||$. In particular, we have
\begin{equation}\label{equivalence}
C_1\|v\|_{H^s(\omgsc)}\leq |||v||| \leq C_2 \|v\|_{H^s(\omgsc)} \;\;\forall\,v\in V_c(\omgsc)
\end{equation}
for some positive constants $C_1$ and $C_2$. As a consequence, any result obtained below involving the energy norm $|||\cdot|||$ can be reinterpreted as a result involving the norm $\|\cdot\|_{H^s(\omgsc)}$. The energy space associated with kernels satisfying \eqref{gamma-conds} and \eqref{case3} is not equivalent to any Sobolev space; however, it is a separable Hilbert space and is a subset of $L^2(\omgsc)$. In both cases, the energy norm satisfies the nonlocal Poincar\'e inequality
\begin{equation}
\|v\|_{L^2(\omgsc)}\leq C_p \, |||v||| \quad \forall\,v\in V_c(\omgsc)
\end{equation}
where $C_p$ is the Poincar\'e constant.

We denote by $V_c^\prime(\omg)$ the dual space of $V_c(\omgsc)$ with respect to the standard $L^2(\omg)$ duality pairing; we define the norm on $V_c^\prime(\omg)$ as
$$
  \|f\|_{V_c^\prime(\omg)}:= \sup_{v\in V_c(\omgsc),\,\,v\ne0}\,\,\frac {\int_\omg fv\,d\xb }{ |||v|||_{V_c(\omgsc)}}.
$$
Note that for Cases 1 and 2 we have that $V_c^\prime(\omg)\subseteq L^2(\omg)$ so that $\|f\|_{V_c^\prime(\omg)}\le \|f\|_{L^2(\omg)}$. In particular, for Case 1 we have that $V_c^\prime(\omg)$ is equivalent to $H^{-s}(\omg)$. We also define the volume ``trace'' space $\widetilde V(\omgc)  := \left\{v|_\omgc : v\in V(\omgsc)\right\}$ and an associated norm 
\begin{equation}\label{tracenorm}
  \|g\|_{\widetilde V(\omgc)} := \inf_{v\in V(\omgsc),\,\,v|_\omgc=g} |||v|||.
\end{equation}

\subsection{Nonlocal volume constrained problems}
We consider the nonlocal volume constrained problem
\begin{equation} \label{eq:fwd_diffusion}
\left\{
\begin{array}{ll}
	\mcD (\vt \mcG u) 	 = f		&	\qquad	\forall \, \xb\in  \omgs\\[2mm]	
	u					 = g 	&	\qquad	\forall \, \xb\in  \omgc
\end{array}\right.
\end{equation}
where we assume that $u\in V_c(\omgsc)$, $f\in V'_c(\omgsc)$, $g\in \widetilde{V}(\omgc)$, and $\vartheta\in \mcC$, where\footnote{The reason of this choice will be made clear in the following section.}
\begin{equation}
\hspace{-2cm}\mcC :=\{\vartheta\in W^{1,\infty}(\omgsc\times\omgsc):\;\; 0<\vartheta_0\leq\vartheta\leq\vartheta_1<\infty, 
         \; \|\vartheta\|_{1,\infty}\leq C<\infty\}.
\end{equation}
Using the nonlocal vector calculus we can define the weak formulation of the problem \eqref{eq:fwd_diffusion} as
\begin{equation}\label{weakf}
\begin{array}{ll}
   &\mbox{\em given $f\in V_c^\prime(\omg)$, $g\in \widetilde V(\omgc)$, and $\vt\in\mcC$, seek $u\in V(\omgsc)$ }\\[1mm]
   &\mbox{\em such that $u=g$ for $\xb\in\omgc$ and}                                                              \\[2mm]
   &\displaystyle\int_\omgsc\int_\omgsc \vt \;\mcG u\cdot\mcG v\,d\yb d\xb = \int_\omg fv\,d\xb
   \qquad\forall\,v\in V_c(\omgsc).
\end{array}
\end{equation}
The problem \eqref{weakf} has a unique solution that depends continuously on the data \cite{dglz:11}.
\section{The optimal control problem}\label{sec:control}
In this section we define the identification problem as a control problem for the nonlocal diffusion equation \eqref{eq:fwd_diffusion}. In a way similar to the local counterpart we demonstrate the existence of at least one solution\footnote{The non-uniqueness of the solution is not due to the nonlocality, the same result holds for the corresponding local PDE control problem.}.

The {\it state} and {\it control} variables are $u$ and the diffusivity $\vt$, respectively, which are related by the {\it state equation} \eqref{eq:fwd_diffusion}; the goal of the control problem is to minimize a cost functional which depends on the state and the control subject to the state equation being satisfied. We define the cost functional as
\begin{equation}\label{eq:cost_func}
J(u,\vartheta):=\displaystyle\frac{1}{2}\int_{\omgs} \big(u-\widehat u \big)^2\,d\xb,
\end{equation}
where $\widehat{u}\in L^2(\omg)$ is a given function. Thus, we want to match as well as possible, in an $L^2(\omg)$ sense, the target function $\widehat{u}$. Formally, we define the control problem as 
\begin{equation}\label{eq:min}
\hspace{-2cm}\begin{minipage}{4.8in}
{\em given $g\in \widetilde V(\omgc)$, $f\in V_c'(\omgsc)$, and $\widehat u\in L^2(\omg)$, seek an optimal control $\vt^\ast\in \mcC$ and an optimal state $u^\ast\in V(\omgsc)$ such that $J(u,\vt)$ given by \eqref{eq:cost_func} is minimized, subject to $u$ and $\vt$ satisfying \eqref{weakf}.} 
\end{minipage}
\end{equation}
\subsection{Existence of an optimal control}
We show that the optimization problem \eqref{eq:min} has at least one solution in the set of admissible parameters; {\color{blue}in this context, the nonlocal Poincar\'e inequality plays a fundamental role.} For simplicity, we consider $g=0$.

{\color{blue}In the proof of the main result we use the following lemma.
\begin{lemma}\label{infty-continuity}
There exists a positive constant $C$ that depends on $\vt_0$, and $f$, such that, for any solutions $u(\vt_i)$ to
\begin{equation}
\left\{\begin{array}{ll}
     \mcD(\vt_i\mcD^*u) = f & \;\;\forall\; \xb\in\omgs \\[1mm]
     u=0                    & \;\;\forall\; \xb\in\omgc
\end{array}\right.
\end{equation}
for $i=a,b$, the following estimate holds:
\begin{equation}\label{nl-homo}
|||u(\vt_a)-u(\vt_b)|||\leq C \|\vt_a - \vt_b\|_\infty.
\end{equation}
\end{lemma}
\noindent{\it Proof.} We first note that 
\begin{displaymath}
\mcD\big(\vt_a\mcD^*(u(\vt_a)-u(\vt_b))\big) = \mcD\big((\vt_b-\vt_a)\mcD^*u(\vt_b)\big);
\end{displaymath}
multiplying both sides by $(u(\vt_a)-u(\vt_b))$ and integrating over $\omgs$ we have
\begin{displaymath}
  \displaystyle\int_\omgs \mcD\big(\vt_a\mcD^*(u(\vt_a)-u(\vt_b))\big) (u(\vt_a)-u(\vt_b)) \,d\xb
= \int_\omgs \mcD\big((\vt_b-\vt_a)\mcD^*u(\vt_b)\big) (u(\vt_a)-u(\vt_b)) \,d\xb.
\end{displaymath}
The nonlocal Green's identity \cite{dglz:11} implies 
\begin{displaymath}
  \displaystyle\int_\omgsc\int_\omgsc \vt_a\big(\mcD^*(u(\vt_a)-u(\vt_b))\big)^2 \,d\yb\,d\xb
= \int_\omgsc\int_\omgsc (\vt_b-\vt_a) \mcD^*u(\vt_b)\cdot \mcD^*(u(\vt_a)-u(\vt_b)) \,d\yb\,d\xb.
\end{displaymath}
Thus,
\begin{displaymath}
\vt_0 |||u(\vt_a)-u(\vt_b) |||^2 \leq \|\vt_a-\vt_b\|_\infty |||u(\vt_b) |||  \; |||u(\vt_a)-u(\vt_b) |||;
\end{displaymath}
then, because the solution of \eqref{nl-homo} depends continuously upon the data, dividing both sides by $|||u(\vt_a)-u(\vt_b) |||$ we have
\begin{displaymath}
|||u(\vt_a)-u(\vt_b) ||| \leq \displaystyle\frac{\|f\|_{V'_c(\omgs)}}{\vt_0}\|\vt_a-\vt_b\|_\infty.
\end{displaymath}
$\boxempty$}
\begin{thm}
There exist at least one solution of the optimization problem \eqref{eq:min}.
\end{thm}
\noindent{\it Proof.} The steps of the proof have been inspired by \cite{Banks} (chapter 6) and \cite{Jin}. We drop the explicit reference to the domain $\omgsc$ and denote $\int_{\omgsc}\int_{\omgsc}$ by $\int\int_{\omgsc}$.

We note that $\inf_{\vt\in\mcC} J(u(\vt),\vartheta)$ is bounded over the set $\mcC$ and, thus, there exists a minimizing sequence $\{\vartheta^n\}\subset \mcC$ such that
$$
J(u(\vartheta^n),\vartheta^n)\rightarrow\inf\limits_{\vartheta\in\mcC} J(u(\vt),\vt) \quad {\rm in} \; \mathbb{R},
$$
where $u(\vt^n)$ denotes the solution of \eqref{weakf} with $\vt=\vt^n$. Because $W^{1,\infty}$ is compactly embedded in $L^\infty$, there exist $\vta\in \mcC$ and a subsequence, which we still denote by $\{\vartheta^n\}$, such that (see, e.g., \cite{Rudin})
\begin{equation}\label{eq:Lstrong}
\vartheta^n \rightarrow\vta \quad {\rm in} \; L^\infty.
\end{equation}
Now, let $u ^n:=u (\vartheta^n)$; then, by definition 
\begin{equation}\label{eq:thetan_qn}
\int\int_{\omgsc} \vartheta^n \mcG u ^n\cdot\mcG v \,d\yb d\xb = \int_{\omgs} f\,v \, d\xb \quad \forall \,v \in V_c.
\end{equation}
The well-posedness of problem \eqref{weakf} implies that there exists $\ua\in V_c$ and a subsequence, which we still denote by $\{u ^n\}$, such that 
\begin{equation}\label{eq:weak}
u ^n\xrightarrow{w}\ua \quad {\rm in} \;V_c,
\end{equation} 
where by $\xrightarrow{w}$ we mean weak convergence (see, e.g., \cite{Rudin}). Now we show that $\ua = u (\vta)$. For all $v\in V_c$ and for all $n$, \eqref{eq:thetan_qn} is equivalent to 
\begin{equation}\label{eq:thetan_qn_equivalent}
\hspace{-2cm}\begin{array}{l}
\displaystyle\int\int_{\omgsc} (\vartheta^n-\vta) \mcG u ^n\cdot\mcG v \,d\yb d\xb + 
\displaystyle\int\int_{\omgsc} \vta \mcG (u ^n-\ua)\cdot\mcG v \,d\yb d\xb \\[3mm]
+\displaystyle\int\int_{\omgsc} \vta \mcG \ua\cdot \mcG v \,d\yb d\xb = 
\int_{\omgs} f\,v \, d\xb.
\end{array}
\end{equation}
Property \eqref{eq:weak} implies that
\begin{equation}
\int\int_{\omgsc} \vta \mcG (u ^n-\ua)\cdot\mcG v \,d\yb d\xb \rightarrow 0
\end{equation}
as $n\rightarrow \infty$. Furthermore, \eqref{eq:Lstrong} implies that
\begin{equation}
\displaystyle\left|\int\int_{\omgsc} (\vartheta^n-\vta) \mcG u ^n\cdot\mcG v \,d\yb d\xb\right| \leq
\displaystyle\|\vartheta^n-\vta\|_\infty \int\int_{\omgsc} \left|\mcG u ^n\cdot\mcG v\right| \,d\yb d\xb
\;\rightarrow 0
\end{equation}
as $n\rightarrow\infty$.
Hence, taking the limit in \eqref{eq:thetan_qn_equivalent} as $n\rightarrow\infty$ we have
\begin{equation}
\int\int_{\omgsc} \vta \mcG \ua \cdot\mcG v \,d\yb d\xb = \int_{\omgs} f\,v \, d\xb
\end{equation}
that gives, by definition $\ua = u(\vta)$. {\color{blue}Next, we note that by choosing $\vt_a=\vt^n$ and $\vt_b=\vta$ in Lemma \ref{infty-continuity} we have that $\vt^n\to\vta$ in $L^\infty$ as $n\to\infty$ implies $u^n\to\ua$ in $V_c$; furthermore, the nonlocal Poincar\'e inequality implies $u^n\to\ua$ in $L^2$ as $n\to\infty$.} Thus,
\begin{equation}\label{eq:proof}
\begin{array}{ll}
J(\ua,\vta) & = \displaystyle\frac{1}{2} \int_{\omgsc} \left(\ua-\widehat u \right)^2\,d\xb = 
                \lim\limits_{n\rightarrow\infty}\frac{1}{2} \int_{\omgsc} \left(u^n-\widehat u \right)^2\,d\xb \\[4mm]
            & = \lim\limits_{n\rightarrow\infty} J(u^n,\vartheta^n) = \inf\limits_{\vartheta\in\mcC} J(u(\vt),\vartheta).
\end{array}
\end{equation} $\boxempty$

Next, for the optimal control problem \eqref{eq:min}, we give a necessary condition for the optimality of a solution. The Lagrangian functional for the problem \eqref{eq:min} is defined as
\begin{equation}\label{eq:lagr}
\hspace{-2cm}\begin{array}{l}
	\displaystyle L(u ,w,\vt,\mu_0,\mu_1) = J(u,\vt) + \int_{\omgs} \big(-\mcD (\vartheta\mcG u ) + f \big) w\,d\xb  +\\
   	                          \hspace{2cm}\displaystyle \int_{\omgsc}\int_{\omgsc}(\vt-\vt_0)\mu_0 \,d\yb d\xb + 
   	                            \int_{\omgsc}\int_{\omgsc}(\vt_1-\vt)\mu_1 \,d\yb d\xb,
\end{array}
\end{equation}
where the {\it adjoint} variable $w\in V_c(\omgsc)$ and $\mu_0,\,\mu_1\in \mcC$ are the Lagrangian multipliers. If $(\ua,\vta)$ is an optimal solution of \eqref{eq:min}, then, $\ua,\,w^*,\,\vta,\,\mu_0^*,$ and $\mu_1^*$ satisfy, respectively, the state, adjoint, optimality, and complementary equations \cite{itku:08}
\begin{subequations}
\begin{equation}
\int_{\omgsc}\int_{\omgsc} \vartheta \mcG \ua \cdot\mcG v \,d\yb d\xb - \int_{\omgs} f\,v \, d\xb=0 \;\;\;\; \forall \,v \in V_c(\omgsc)\label{eq:weak_state}\qquad\qquad\qquad\qquad\\
\end{equation}
\begin{equation}
\int_{\omgsc}\int_{\omgsc} \vartheta \mcG w^* \cdot\mcG\psi \,d\yb d\xb - \int_{\omgs} (\ua-\widehat u)\,\psi \, d\xb=0 \;\;\;\; \forall \,\psi \in V_c(\omgsc)\label{eq:weak_adjoint}\qquad\qquad\\
\end{equation}
\begin{equation}
\int_{\omgsc}\int_{\omgsc}  \left(\mcG \ua  \cdot\mcG w^*+ \mu_0^* - \mu_1^*\right)\varphi \,d\yb d\xb =0 \;\;\;\; \nonumber\\ \forall \varphi \in W^{1,\infty}(\omgsc\times\omgsc) \label{eq:weak_opt}\\
\end{equation}
\begin{equation}
\int_{\omgsc}\int_{\omgsc}  \mu_0^*\left(\vt^* - \vt_0\right) \,d\yb d\xb = 0 \label{eq:compl1}\\
 \end{equation}
\begin{equation}
\int_{\omgsc}\int_{\omgsc}  \mu_1^*\left(\vt_1 - \vt^*\right) \,d\yb d\xb = 0. \label{eq:compl}
\end{equation}
\end{subequations}

\section{Finite dimensional approximation}\label{fin_dim_approx}
In this section we consider the convergence of solutions of finite dimensional discretizations of the optimal control problem. For finite element discretizations we also derive a priori error estimates for the state and control variables. We limit ourselves to Case 1 and we analyze the homogeneous Dirichlet problem.

We choose the families of finite dimensional subspaces
\begin{equation}\label{conf}
V^N(\omgsc) \subset V(\omgsc), \qquad W^M(\omgsc)\subset W^{1,\infty}(\omgsc\times\omgsc)
\end{equation}
parametrized by integers $N,\,M\to\infty$, and then define the constrained finite dimensional subspace
\begin{equation}\label{confc}
V^N_c(\omgsc) := V^N(\omgsc)\cap V_c(\omgsc)=
\left\{v\in V^N(\omgsc) \,\,:\,\, v=0\;{\rm on}\;\omgc\right\}.
\end{equation}
The usual choice for $N$ and $M$ is the dimension of the subspaces. We assume that, for any function $v\in V(\omgsc)$ and any function $\sigma\in W^{1,\infty}(\omgsc\times\omgsc)$, the sequence of best approximations with respect to the energy norm $|||\cdot|||$ and the $W^{1,\infty}$ norm, respectively, converges, i.e., 
\begin{equation}\label{baerror}
  \lim\limits_{N\to\infty} \,  \inf\limits_{v_N\in V^N} ||| v - v_N||| = 0
  \qquad\forall\, v \in V(\omgsc)
\end{equation}
and
\begin{equation}\label{baerrorW}
\lim\limits_{M\rightarrow\infty}  \, \inf\limits_{\sigma_M\in W^M}\|\sigma_M-\sigma\|_{1,\infty} = 0
\qquad\forall\, \sigma \in W^{1,\infty}(\omgsc\times\omgsc).
\end{equation}
The admissible parameter sets are specified as
\begin{equation}
\hspace{-2cm}\mcC^M=\{\vt_M\in W^M(\omgsc)\; {\rm s.t.} \;\; 0<\vt_0\leq\vt_M\leq\vt_1<\infty, \;\; \|\vt_M\|_{1,\infty}<C<\infty\}.
\end{equation} 
We seek the Ritz-Galerkin approximations $u_N \in V^N_c(\omgsc)$ and $\vt_M\in\mcC^M$ determined by posing \eqref{eq:min} on $V^N_c(\omgsc)$ and $W^M(\omgsc)$. The finite dimensional state equation in a weak form is given by
\begin{equation}\label{eq:weak_finite}
\hspace{-2cm}\int_{\omgsc}\int_{\omgsc} \vtn \mcG u_N \cdot\mcG v_N \,d\yb d\xb = \int_{\omgs} f\,v_N \, d\xb \qquad \forall \,v_N \in V^N_c(\omgsc).
\end{equation}
As for the infinite dimensional case, by the Lax-Milgram theorem, for all $N$ and $M$, \eqref{eq:weak_finite} has a unique solution $u_N \in V_c^N(\omgsc)$ for all $\vtn\in\mcC^M$. Then, the finite dimensional control problem is formulated as
\begin{equation} \label{eq:min_finite}
\hspace{-2cm}	\min\limits_{\vtn\in\mcC^M}	 J(u_N,\vtn)	 :=\displaystyle\frac{1}{2}\int_{\omgs} \big(u_N - \widehat u \big)^2\,d\xb
	\quad\mbox{such that \eqref{eq:weak_finite} is satisfied}.
\end{equation}
Using the same arguments as for the infinite dimensional problem it is possible to show that problem \eqref{eq:min_finite} has at least one solution. 

Next, we consider a finite element approximation for the case that both $\omgsc$ and $\omgs$ are polyhedral domains. We partition $\omgsc$ and $\omgsc\times\omgsc$ into finite elements for the discretization of $u$ and $\vt$ respectively and we denote by $h_u$ and $h_\vt$ the diameter of the largest element in each partition. We assume that the partitions are shape-regular and quasiuniform \cite{brsc:08} as the grid sizes $h_u,\,h_\vt\to 0$, i.e., as $N,\,M\to\infty$. We choose $V_{c}^N(\omgsc)$ and $W^M(\omgsc)$ to consist of piecewise polynomials of degree no more than $m$ and $l$, respectively, defined with respect to each grid.  For some real constants $K_1,\,K_2$, and $K_3$, the following assumptions are made for the finite dimensional subspaces $V_c^N(\omgsc)$ and $W^M(\omgsc)$. For all $v\in V_c(\omgsc)\cap H^{m+t}(\omgsc)$,  $s\in (0,\,1)$ and $t\in[s,\,1]$, there exists $\Pi^Vv\in V^N_c(\omgsc)$ such that
\begin{equation}\label{eq:piV}
\|v-\Pi^Vv\|_{H^s(\omgsc)}\leq K_1 N^{-(m+t-s)}\|v\|_{H^{m+t}(\omgsc)}.\\
\end{equation}
For all $\sigma\in W^{l+1,\infty}(\omgsc\times\omgsc)$ there exists $\Pi^W\sigma\in W^M(\omgsc)$ such that
\begin{equation}\label{eq:piW}
\|\sigma-\Pi^W\sigma\|_{L^{\infty}(\omgsc\times\omgsc)}\leq K_2 M^{-(l+1)}\|\sigma\|_{l+1,\infty(\omgsc)}.\\
\end{equation}
Also, we assume that for all $v_N\in V_c^N(\omgsc)$ the inverse inequality  
\begin{equation}\label{eq:inverse_ineq}
\|v_N\|_{H^s(\omgsc)}\leq K_3\,N^s\,\|v_N\|_{L^2(\omgsc)}
\end{equation}
holds. These properties are satisfied for  wide classes of finite element spaces; see \cite[p. 121]{ciarlet} for \eqref{eq:piW} and \cite{ciarlet, graham} for \eqref{eq:inverse_ineq}.
Throughout this section we let $K$ denote a generic constant, independent on $N$ and $M$, and we suppress explicit reference to the domain $\omgsc$. The proofs of lemmas and theorems are the nonlocal equivalent of \cite[Theorem VI.3.1]{Banks}. 

First, we consider the following two lemmas.
\begin{lemma}\label{lemma1}
Let $\vta\in\mcC$ be a solution of \eqref{eq:min} and let $\ua$ be the corresponding solution of \eqref{weakf}. Then,
\begin{equation}\label{eq:lemma1}
\hspace{-2cm}|||u(\Pi^W \vta)-\Pi^V\ua |||\leq \frac{C_2 C_u}{\vt_0} \|\Pi^W\vta-\vta\|_{L^\infty(\omgsc)} + \frac{\vt_1}{\vt_0} |||\Pi^V\ua-\ua |||,
\end{equation} 
where $C_2$ is the equivalence constant in \eqref{equivalence} and $C_u$ is such that $\|u^*\|_{H^s(\omgsc)}\leq C_u$.
\end{lemma}
{\it Proof.}
Let $\widetilde\vt_M=\Pi^W\vta$ and $\unt=u(\vtnt)$. Observe that
\begin{equation}
\int\int_{\omgsc}\vtnt\mcG \unt\cdot\mcG v_N \,d\yb d\xb= 
\int_{\omgs} f\,v_N \,d\xb=\displaystyle\int\int_{\omgsc}\vta\mcG \ua\cdot\mcG v_N \,d\yb d\xb 
\end{equation}
for all $v_N\in V_c^N$. Then,
\begin{equation}
\begin{array}{ll}
   & \displaystyle\int\int_{\omgsc}\vtnt\mcG (\unt-\Pi^V\ua)\cdot\mcG v_N \,d\yb d\xb \\[3mm]
 = & \displaystyle\int_{\omgs} f\,v_N \,d\xb- \int\int_{\omgsc}\vtnt\mcG (\Pi^V\ua)\cdot\mcG v_N \,d\yb d\xb \\[3mm]
 = & \displaystyle\int\int_{\omgsc}\vta\mcG \ua\cdot\mcG v_N \,d\yb d\xb - 
     \int\int_{\omgsc}\vtnt\mcG (\Pi^V\ua)\cdot\mcG v_N \,d\yb d\xb \\[3mm]
 = & \displaystyle\int\int_{\omgsc}(\vta-\vtnt)\mcG \ua\cdot\mcG v_N \,d\yb d\xb -  
     \int\int_{\omgsc}\vtnt\mcG (\ua-\Pi^V\ua)\cdot\mcG v_N \,d\yb d\xb. 
\end{array}
\end{equation}
Now, we choose $v_N=\unt-\Pi^V\ua$, thus
\begin{equation}
\begin{array}{ll}
        & \vt_0\,|||\unt-\Pi^V\ua |||^2                                                                              \\[2mm]
   \leq & \left|\displaystyle\int\int_{\omgsc}\vtnt\mcG (\unt-\Pi^V\ua)\cdot\mcG (\unt-\Pi^V\ua) \,d\yb d\xb \right| \\[3mm]
   \leq & \displaystyle\left|\int\int_{\omgsc}(\vta-\vtnt)\mcG \ua\cdot\mcG (\unt-\Pi^V\ua) \,d\yb d\xb \right|      \\[3mm]
   +    & \displaystyle\left|\int\int_{\omgsc}\vtnt\mcG (\ua-\Pi^V\ua)\cdot\mcG (\unt-\Pi^V\ua) \,d\yb d\xb\right|   \\[3mm]
   \leq & \displaystyle\|\vtnt-\vta\|_{L^\infty} \,|||\ua |||\,|||\unt-\Pi^V\ua ||| + \|\vtnt\|_{L^\infty}
          \,|||\ua-\Pi^V\ua |||\,|||\unt-\Pi^V\ua |||                                                                \\[2mm]
   \leq & \displaystyle C_2 C_u \|\vtnt-\vta\|_{L^\infty} |||\unt-\Pi^V\ua ||| + \vt_1 |||\ua-\Pi^V\ua |||\,|||\unt-\Pi^V\ua |||.   
\end{array}
\end{equation}
Dividing both sides by $|||\unt-\Pi^V\ua |||$ we obtain \eqref{eq:lemma1}. $\boxempty$
\begin{lemma}\label{lemma2}
Let $m$ and $l$ be non-negative integers and let $s\in(0,\,1)$ and $t\in[s,\,1]$. Then, there exists a constant $K$, independent on $N$ and $M$, such that, for sufficiently large $N$ and $M$ and for sufficiently smooth $u^*$ and $\vt^*$ (defined as in Lemma \ref{lemma1})
\begin{equation}\label{lemma2res}
\hspace{-2cm}\inf\limits_{\vt_M\in W^M} \|\ua-u_N(\vtn)\|_{L^2} \leq 
K \left(M^{-(l+1)}\|\vta\|_{l+1,\infty}+ N^{-(m+t-s)} \|\ua\|_{H^{m+t}}\right).
\end{equation} 
\end{lemma}
\noindent{\it Proof.}
Let $\widehat u =\unt$ (defined as in Lemma \ref{lemma1}); then the optimal adjoint variable $w^*$, the solution of \eqref{eq:weak_adjoint}, satisfies \cite{dglz:11}
\begin{equation}\label{eq:fwd_ineq}
\|w^*\|_{H^s}\leq K \|\unt-\ua\|_{L^2}.
\end{equation}
We consider equation \eqref{eq:weak_adjoint} and choose the test function $\unt-\ua$; we then have
$$
\begin{array}{ll}
       & \|\unt-\ua\|^2_{L^2}= \displaystyle\int_{\omgs} (\unt-\ua)(\unt-\ua)\,d\xb                                          \\ [3mm]
  =    & \displaystyle \int\int_{\omgsc}\vta \mcG w^*\cdot\mcG (\unt-\ua)\,d\yb d\xb                                         \\[3mm]
  =    & \displaystyle\int\int_{\omgsc}(\vta-\vtnt) \mcG w^*\cdot\mcG (\unt-\ua)\,d\yb d\xb                                  \\[3mm]
  +    & \displaystyle\int\int_{\omgsc}\vtnt \mcG w^*\cdot\mcG (\unt-\ua)\,d\yb d\xb                                         \\[3mm]
  \leq & \displaystyle\|\vta-\vtnt\|_{L^\infty}\,|||w^*|||\,|||\unt-\ua||| + \|\vtnt\|_{L^\infty}\,|||w^*|||\,|||\unt-\ua||| \\[2mm]
  \leq & \displaystyle K \left(\|\vta-\vtnt\|_{L^\infty}\|\unt-\ua\|_{L^2}|||\unt-\ua||| 
  +      \vt_1\,\|\unt-\ua\|_{L^2}|||\unt-\ua|||\right),
\end{array}
$$
where we used \eqref{eq:fwd_ineq}. Dividing both sides by $\|\unt-\ua\|_{L^2}$, we have
$$
\begin{array}{ll}
    \|\unt-\ua\|_{L^2} & \leq K |||\unt-\ua ||| \left(\|\vta-\vtnt\|_{L^\infty} + \,\vt_1 \right)\\[2mm]
    & \leq K\left(|||\unt - \Pi^V\ua||| + ||| \Pi^V\ua-\ua||| \right) \left(\|\vta-\vtnt\|_{L^\infty} + \,\vt_1 \right).
\end{array}
$$
Using Lemma \ref{lemma1}, \eqref{eq:piV}, and \eqref{eq:piW}, we obtain
$$
\begin{array}{l}
     \|\unt-\ua\|_{L^2}\leq \left(\displaystyle\frac{C_2 C_u}{\vt_0} \|\vtnt-\vta\|_{L^\infty} 
                          + \left(1+\frac{\vt_1}{\vt_0}\right) |||\Pi^V\ua-\ua ||| \right)              \\[2mm]
     \qquad\left(\|\vta-\vtnt\|_{L^\infty} + \,\vt_1 \right)                                            \\[5mm]
     \leq K\left( \|\vtnt-\vta\|^2_{L^\infty} + \|\vtnt-\vta\|_{L^\infty}\right.                        \\[1mm]
     \qquad\left. + \|\vtnt-\vta\|_{L^\infty} |||\Pi^V\ua-\ua ||| +  |||\Pi^V\ua-\ua ||| \right)        \\[5mm]
     \leq K \left(K^2_2 M^{-2(l+1)}\|\vta\|^2_{l+1,\infty}+ K_2 M^{-(l+1)}\|\vta\|_{l+1,\infty} \right. \\[2mm]
     \qquad\left. + K_1 K_2 M^{-(l+1)}N^{-(m+t-s)}\|\vta\|_{l+1,\infty} \|\ua\|_{H^{m+t}} + 
                    K_1 N^{-(m+t-s)}\|\ua\|_{H^{m+t}} \right)                                           \\[5mm]
     \leq K \left(M^{-(l+1)}\|\vta\|_{l+1,\infty}+ N^{-(m+t-s)} \|\ua\|_{H^{m+t}}\right).
\end{array}
$$
Now, because $\inf_{\vt_M\in W^M}  \|\ua-u_N(\vtn)\|_{L^2}\leq \|\unt-\ua\|_{L^2}$, \eqref{lemma2res} follows. $\boxempty$

Finally, we can state the main theorem, which provides an estimate of the approximation error for the control variable $\vt$.
\begin{thm}\label{th:param_apriori}
Let $m$ and $l$ be non-negative integers and let $s\in(0,\,1)$ and $t\in[s,\,1]$. Assume that for $f\in H^{m+t}(\omgs)$ and $\widehat u \in L^2(\omgs)$, $\vta\in W^{l+1,\infty}(\omgsc\times\omgsc)$ is a solution of \eqref{eq:min} and $u(\vta) = \ua \in V_c(\omgsc) \cap H^{m+t}(\omgsc)$ is the corresponding state. Then, there exists a constant $K$, independent on $N$ and $M$ such that, for every solution $\vtnb$ of \eqref{eq:min_finite} 
\begin{equation}\label{main_thm}
\begin{array}{ll}
   \qquad\displaystyle\int\int_{\omgsc}\left|\vta-\vtnb\right|\mcG \ua\cdot\mcG \ua\,d\yb d\xb & 
   \leq \;K \left(N^s dist(\widehat u ,\mcV) \right.\\[2mm]
   & +\; N^{-(m+t-2s)} \|\ua\|_{H^{m+t}(\omgsc)}\\[2mm]
   & +\;N^s M^{-(l+1)} \|\vta\|_{l+1,\infty} \\[2mm]
   & + \left. N^{-(m+t-s)} \|f\|_{H^{m+t}(\omgsc)}\right),
\end{array}
\end{equation}
where $\mcV = \{u(\vt)\,:\, \vt\in\mcC\}$. Furthermore, if $\widehat u \in \mcV$, then
\begin{equation}
\begin{array}{ll} 
     \displaystyle\int\int_{\omgsc}\left|\vta-\vtnb\right|\mcG \ua\cdot\mcG \ua\,d\yb d\xb
     & \leq \;K \left(N^{-(m+t-2s)} \|\ua\|_{H^{m+t}(\omgsc)} \right.                         \\[2mm]
     & + \;N^s M^{-(l+1)} \|\vta\|_{l+1,\infty}                                               \\[2mm]
     & + \left. N^{-(m+t-s)} \|f\|_{H^{m+t}(\omgsc)}\right).
\end{array}
\end{equation}
\end{thm}
\noindent{\it Proof.}
Let  $R_1 := \{(\xb,\,\yb) : \vta - \vtnb \geq 0 \}$ and $R_2 := (\omgsc\times\omgsc)\backslash R_1$ and define the function $\chi:\omgsc\times\omgsc\rightarrow\mathbb{R}$ by
$$
\chi(\xb,\yb) := \left\{
\begin{array}{rl}
1 & (\xb,\,\yb)\in R_1\\
-1 & (\xb,\,\yb)\in R_2.
\end{array}\right.
$$  
Recall that 
$$
\mcD(\vta \mcG\ua) = f \quad \hbox{and} \quad \mcD(\vtnb \mcG u^*_N) = \Pi^N f,
$$
where by $\Pi^N$ we denote the orthogonal $L^2$-projection onto $V_c^N(\omgsc)$. Therefore
$$
\mcD((\vta-\vtnb) \mcG\ua) = \mcD(\vtnb \mcG (u_N^*-\ua)) + f-\Pi^N f.
$$
Taking the inner product with the function $\ua\chi$, we have
$$
\begin{array}{l}
      \displaystyle-\int_{\omgs}\mcD(|\vta-\vtnb| \mcG \ua)\ua\,d\yb d\xb\\[3mm]
   =  \displaystyle \int_{\omgs}\mcD(\vtnb \mcG(u_N^*-\ua))\ua\chi \,d\yb d\xb + \int_{\omgs}(f-\Pi^N f)\ua\chi\,d\xb.
\end{array}
$$
Thus, we write
\begin{equation}\label{eq:first_ineq}
\begin{array}{ll}
        & \displaystyle\int\int_{\omgsc}|\vta-\vtnb|\mcG \ua \cdot\mcG \ua\,d\yb d\xb                          \\[3mm]
   \leq &\left|\displaystyle\int\int_{\omgsc}\vtnb \mcG (u_N^*-\ua)\cdot\mcG (\ua\chi) \,d\yb d\xb\right| 
        + \left|\displaystyle\int_{\omgs}(f-\Pi^N f)\ua\chi\,d\xb\right|                                       \\[3mm]
   \leq & \|\vtnb\|_{L^{\infty}}|||u_N^*-\ua |||\,|||\ua\chi||| + \|f-\Pi^N f\|_{L^2}\|\ua\chi\|_{L^2}         \\[2mm]
   \leq & \vt_1 C_2\|\ua\|_{H^s}  |||u_N^*-\ua||| + \|f-\Pi^N f\|_{L^2}\|\ua\|_{H^s}                           \\[2mm]
   \leq & \vt_1C_2 C_u  |||u_N^*-\ua||| + C_u \|f-\Pi^N f\|_{H^s}                                              \\[2mm]
   \leq & \vt_1 C^2_2 C_u \|u_N^*-\ua\|_{H^s} + C_u K_1 N^{-(m+t-s)} \|f\|_{H^{m+t}},
\end{array}
\end{equation} 
where we used \eqref{eq:piV}. Next, we find a bound for $\|u_N^*-\ua\|_{H^s}$:
$$
\begin{array}{ll}
     \|u_N^*-\ua\|_{H^s} & \leq \|\ua-\Pi^V\ua\|_{H^s} +\|u_N^*-\Pi^V\ua\|_{H^s}                                 \\[2mm]
   & \leq K_1 N^{-m-t+s} \|\ua\|_{H^{m+t}} + K_3 N^s \|u_N^*-\Pi^V\ua\|_{L^2}                                    \\[2mm]
   & \leq K_1 N^{-m-t+s} \|\ua\|_{H^{m+t}} + K_3 N^s \left(\|\ua-\Pi^V\ua\|_{L^2} + \|\ua-u_N^*\|_{L^2}\right)   \\[2mm]
   & \leq \left(K_1 N^{-m-t+s} + K_1 K_3 N^{-(m+t-2s)}\right) \|\ua\|_{H^{m+t}} + K_3 N^s\|\ua-u_N^*\|_{L^2},
\end{array}
$$
where we used  \eqref{eq:piV} and \eqref{eq:inverse_ineq}. Then, we find a bound for $\|\ua-u_N^*\|_{L^2}$:
\begin{equation}\label{eq:l2_error_norm}
\begin{array}{ll}
     \|\ua-u_N^*\|_{L^2} & \leq \|\ua-\widehat u \|_{L^2} +\|\widehat u -u_N^*\|_{L^2}                               \\
   & \leq dist(\widehat u ,\mcV) +\inf\limits_{W^M}\|\widehat u -u_N(\vtn)\|_{L^2}                                   \\
   & \leq dist(\widehat u ,\mcV)+ \inf\limits_{W^M}\left(\|\widehat u -\ua\|_{L^2}+ \|\ua-u_N(\vtn)\|_{L^2}\right)   \\
   & \leq 2dist(\widehat u ,\mcV) + \inf\limits_{W^M} \|\ua-u_N(\vtn)\|_{L^2}                                        \\
   & \leq 2dist(\widehat u ,\mcV)+ K M^{-(l+1)}\|\vta\|_{l+1,\infty} + K N^{-(m+t-s)} \|\ua\|_{H^{m+t}(\omgsc)},
\end{array}
\end{equation}
where we applied Lemma \ref{lemma2}. Thus,
\begin{equation}\label{eq:Hs_estimate}
\begin{array}{ll}
\|u_N^*-\ua\|_{H^s} & \leq 2\,K_3\,dist(\widehat u ,\mcV)+K N^sM^{-(l+1)}\|\vta\|_{l+1,\infty} \\[2mm]
                    & +\left( K_1 N^{-(m+t-s)}+ K_1 K_3 N^{-(m+t-2s)} +K K_3 N^{-(m+t-2s)}\right)\|\ua\|_{H^{m+t}(\omgsc)}.
\end{array}
\end{equation}
Combining \eqref{eq:Hs_estimate} and \eqref{eq:first_ineq} we obtain \eqref{main_thm}. $\boxempty$

Note that Theorem \ref{th:param_apriori} provides an estimate for the $L^2$ norm of the approximation error for the state variable. In fact, from equation \eqref{eq:l2_error_norm}, we have 
\begin{equation}\label{eq:l2_error}
\hspace{-2cm}\|u_N^*-\ua\|_{L^2}\leq K\left( dist(\widehat u ,\mcV)+ M^{-(l+1)}\|\vta\|_{l+1,\infty} + N^{-(m+t-s)} \|\ua\|_{H^{m+t}(\omgsc)}\right).
\end{equation}

\section{Numerical tests}\label{num_tests}
In this section we present the results of computational experiments for finite element discretizations of one-dimensional problems. These preliminary results illustrate the theoretical results in Section \ref{fin_dim_approx} and provide the basis for extensions to two- and three-dimensional experiments. 

{\color{blue} We conduct two convergence analyses. First, for both Case 1 and Case 2 we consider the convergence of approximate optimal solutions to fine-grid surrogates for the analytic optimal solutions; here, we do not assume any knowledge of the optimal solution in $\omgc$, where we prescribe homogeneous Dirichlet conditions. Then, for Case 2, we analyze the convergence to optimal analytic solutions prescribing exact volume constraints in $\omgc$. We describe the problem settings used in our tests.
\paragraph{Case 1} We consider $\omgs=(-1,1)$, $\omgc=(-1-\varepsilon,-1)\cup(1,1+\varepsilon)$,
\begin{equation}\label{kernel1}
\gamma_1(x,y) = \dfrac{1}{|x-y|^{1+2s}}, \qquad s\in(0,1),
\end{equation}
for which $\mcL$ corresponds to a fractional differential operator \cite{mesi:11}, and the data set
\begin{equation}
{\rm A:}
\left\{\begin{array}{l}
f(x) = 1                       \\[1mm]
\widehat{u}(x) = u_{\rm A}(x)  \\[1mm]
\left.u(x)\right|_{\omgc} = 0, \\[1mm]
s = 0.7.
\end{array}\right.  
\end{equation}
Here $u_{\rm A}(x)$ is a surrogate for an exact solution of \eqref{eq:fwd_diffusion}; it corresponds to the finite element approximation computed on a very fine grid using
\begin{equation}
\vt_{\rm A}(x,y) = 2+0.4(x+y-1)^2.
\end{equation}
\paragraph{Case 2} We consider $\omgs=(0,1)$, $\omgc=(-\varepsilon,0)\cup(1,1+\varepsilon)$, and
\begin{equation}\label{kernel2}
\gamma_2(x,y) = \dfrac{1}{\varepsilon^2|x-y|},
\end{equation}
which is often used in the literature, e.g. in a linearized model for continuum mechanics \cite{chgu:11}. We introduce the following data sets
\begin{equation}
\begin{array}{lll}
{\rm B:}\left\{\begin{array}{l}
    f(x) = \varepsilon^2 + 24 x^2 - 24 x +7.6  \\[1mm]
    \widehat{u}(x) = 2.5\,x(1-x)               \\[1mm]
    \left.u(x)\right|_{\omgc} = 2.5\,x(1-x)
\end{array}\right.
&
\quad{\rm C:}\left\{\begin{array}{l}
    f(x) = 5                         \\[1mm]
    \widehat{u}(x) = u_{\rm C}(x)    \\[1mm]
    \left.u(x)\right|_{\omgc} = 0
\end{array}\right.
& 
\quad{\rm D:}\left\{\begin{array}{l}
    f(x) = 5                         \\[1mm]
    \widehat{u}(x) = u_{\rm D}(x)    \\[1mm]
    \left.u(x)\right|_{\omgc} = 0.
\end{array}\right.
\end{array}
\end{equation}
In case B, $\widehat{u}$ is the solution of \eqref{eq:fwd_diffusion} for $\vt(x,y) = \vt_{\rm A}$; thus, $\widehat{u}$ and $\vt_{\rm A}$ are the optimal state and parameter. In cases C and D, $u_{\rm C}$ and $u_{\rm D}$ are surrogates for an exact solution of \eqref{eq:fwd_diffusion}; they are, in fact, the finite element solutions computed on a very fine grid using respectively $\vt(x,y) = \vt_{\rm C}(\frac{x+y}{2})$ and $\vt(x,y) = \vt_{\rm D}(\frac{x+y}{2})$, where
\begin{equation}\label{eq:piecewiseT}
\vt_{\rm C}(z) = \left\{\begin{array}{ll}
     0.2+(z-0.625)^2      & \; z\in (0,0.625)      \\[1mm]
     z+1.25               & \; z\in (0.625,0.75)   \\[1mm]
     14.4(z-0.75)+2       & \; z\in(0.75,1)
\end{array}\right.\qquad\quad
\vt_C(z) = \left\{\begin{array}{ll}
     1   & \; z\in(0,0.2)     \\[1mm]
     0.1 & \; z\in(0.2,0.6)   \\[1mm]
     1   & \; z\in(0.6,1).
\end{array}\right.
\end{equation}}
\paragraph{Implementation details} 
For the state variable we introduce a partition of $\overline\omgsc = \overline{(a-\varepsilon,  b+ \varepsilon)}$ such that, for the positive constants $K_u$ and $J_u$
\begin{equation}\label{eq:partitionu}
\begin{array}{ll}
    &  a - \varepsilon = x_{-K_u} < \cdots < x_{-1} < a = x_0 < x_1 < \cdots < x_{J_u-1}\\
    &  \qquad\qquad <x_{J_u} = b<x_{J_u+1} < \cdots < x_{J_u+K_u} = b + \varepsilon.
\end{array}
\end{equation} 
Then, $h_u$ is defined as $\max_{j = -K_u,\ldots,K_u+J_u-1}|x_{j+1}-x_j|$. In the numerical experiments we let $V^N(\omgsc)$ be the finite element space of piece-wise linear polynomials. 

\begin{figure}[h!]
\begin{center}
\includegraphics[width=0.5\textwidth]{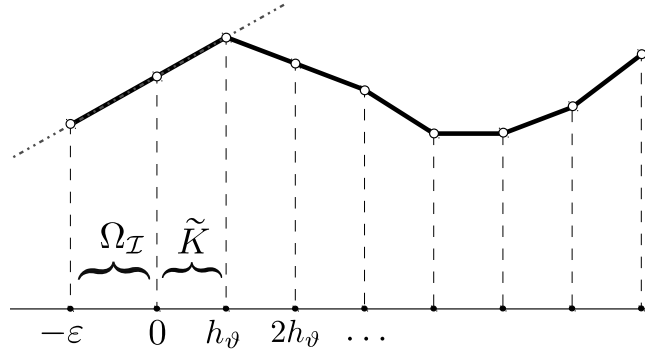} 
\end{center}
\caption{Example of a function in $W^M(\omgsc)$ for the data sets B and C.}
\label{WM_example}
\end{figure}

For the approximation of the parameter we consider $\vt(x,y)$ as a function of one variable only, i.e., as in \eqref{eq:piecewiseT}, $\vt(x,y) = \vt(\frac{x+y}{2})$, defined in $\omgsc$. Thus, for the positive constants $K_\vt$ and $J_\vt$ we introduce the partition
\begin{equation}\label{eq:partitionvt}
\begin{array}{ll}
    &  a - \varepsilon = x_{-K_\vt} < \cdots < x_{-1} < a = x_0 < x_1 < \cdots < x_{J_\vt-1}\\
    &  \qquad\qquad < x_{J_\vt}=1<x_{J_\vt+1} < \cdots < x_{J_\vt+K_\vt} = b + \varepsilon.
\end{array}
\end{equation} 
Then, $h_\vt$ is defined as $\max_{j = -K_\vt,\ldots,K_\vt+J_\vt-1}|x_{j+1}-x_j|$. In choosing $W^M(\omgsc)$ one has to be careful; the most natural choice, for cases A and B and C, is to let $W^M(\omgsc)$ be the space of continuous piece-wise linear polynomials defined over the partition. However, this makes the problem very ill-conditioned. To circumventing this problem we define $W^M(\omgsc)$ as the space of continuous piece-wise linear polynomials such that, for all $\sigma_M\in W^M(\omgsc)$, $\left.\sigma_M\right|_\omgc$ is a linear extension of $\left.\sigma_M\right|_{\widetilde K}$, being $\widetilde K$ the element of the partition adjacent to $\omgc$. An example of a function belonging to $W^M(\omgsc)$ is displayed in Figure \ref{WM_example} for the domain configuration of Case 2.
This empirical choice is not motivated by a theoretical result but by numerical evidence; in fact, we observe that the accuracy of the numerical solutions of the state equation, when the analytic parameter is projected onto $W^M(\omgsc)$, is not affected. 

For case D we define $W^M(\omgsc)$ as the space of piece-wise constant functions such that, for all $\sigma_M\in W^M(\omgsc)$, $\left.\sigma_M\right|_{(a-\varepsilon,a)} = \left.\sigma_M\right|_{(a,a+h_\vt)}$ and $\left.\sigma_M\right|_{(b,b+\varepsilon)} = \left.\sigma_M\right|_{(b-h_\vt,b)}$
The finite dimensional optimization problem is solved with the Broyden{-}Fletcher{-}Goldfarb{-}Shanno (BFGS) method \cite{nocedal99}, without prescribing any conditions on the lower and upper bounds of the parameter. The reason of this choice is because in our simulations those bounds are not violated.

For the solution of local \cite{itku:08,gunz:02} and nonlocal \cite{degu:13} optimization problems a regularization term is usually added to the functional to prevent the ill-posedness and the ill-conditioning of the mathematical and numerical problems. In this case (as in its local counterpart), as shown in Section  \ref{sec:control} and Section \ref{fin_dim_approx}, the problems \eqref{eq:min} and \eqref{eq:min_finite} admit solutions in the space of admissible parameters without regularization. However, a regularization term can be added in case of ill-conditioning; this is the case of the data set D where we utilize the functional
\begin{equation}\label{regularization}
\overline{J}(u_N,\vt_M) = J(u_N,\vt_M) + \beta\sum\limits_{j=1}^{J_\vt -1} \big(\vt_M(\overline{x}_j)-\vt_M(\overline{x}_{j+1})\big)^2, 
\end{equation}
where $\beta>0$, and $\overline{x}_j$ is any point inside the interval $[x_j,x_{j+1})$, for $j=1,\ldots J_\vt -1$. The additional term in \eqref{regularization} has the effect of minimizing the jumps in $\vt_M$.
\subsection{Convergence of the finite element approximate optimal solutions}
{\color{blue} We analyze the convergence with respect to the grid sizes $N$ and $M$ of finite element approximations to the surrogates, for cases A, C, and D, and to the analytic solution for case B.

In Tables \ref{tab:fractional-kernel} and \ref{tab:Hconv_exact} we report the error and the corresponding rates for cases A and B respectively. Here, $e_{u,2} = \|u_N^*-\widehat u\|_{L^2(\omgs)}$ and $e^*_\vt = \|\vt_M^*-\vt_{\rm A}\|_*$, where $\|\cdot\|_*$ is the left-hand side of \eqref{main_thm}; $u_{\rm A}$ is generated using $N=2^{11}$. Note that, because $\widehat{u}\in\mcV$, the term $dist(\widehat{u},\mcV)$ in the estimates \eqref{main_thm} and \eqref{eq:l2_error} is negligible. For case A, though the optimal analytic solution does not necessarily belong to $H^{1+t}(\omgsc)$, $t\in[s,1]$, our results are consistent with the estimate of Theorem \ref{th:param_apriori}; in fact, when the solution is regular enough, for $s=0.7$, we expect to observe a convergence rate $r_u\in(1,1.3)$ for the state and $r_\vt\in(0.3, 0.6)$ for the parameter. 

For case B, where both $u^*$ and $\vt^*$ belong to $C^\infty(\omgsc)$, we observe a quadratic convergence for both the error norms of the state and the parameter for different choices of interaction radius. 
Up to our knowledge, accurate theoretical results concerning the rates of convergence of the discretization error are not available for kernels as in Case 2; however, there is numerical evidence that for piecewise linear finite element approximations the numerical solutions of the state equation do converge to the analytic solutions in $L^2$ as $h\to 0$ with quadratic convergence rate, see the extensive analysis conducted in \cite{chgu:11}. Furthermore, in \cite{degu:13}, numerical tests show that the same convergence rate is preserved for the approximate solutions of optimal control problems constrained by a nonlocal diffusion equation where the control parameter is the source term $f$. 
The plots of the optimal approximate solutions are not significant, because of the superimposition of the solutions, and are not reported.

In case C, $u_{\rm C}$ is determined using $N=2^{12}$. Here, the parameter $\vt_{\rm C}$ has a discontinuous derivative; thus, we do not expect to observe the same convergence rates as in case B. In fact, though we have convergence, the rates do not show a specific trend and for this reason they are not reported. In Figure \ref{fig:tDNM_DG} we report the approximate optimal parameter for several values of $M$ and we observe a very good match with $\vt_{\rm C}$ as we refine the grid. In Figure \ref{fig:uREF_DG} (left) we report $u_{\rm C}$.
\begin{table}[t]
\begin{center}
\begin{tabular}{|c|c|c|c|c|c|}
\multicolumn{6}{c}{$\varepsilon = 2^{-4}$}          \\ \hline
$N$   & $M$   & $e_{u,2}$ & rate & $e^*_\vt$ & rate \\ \hline
$2^4$ & $2^2$ & 2.44e-03  & -    & 2.35e-02  & -    \\ 
$2^5$ & $2^3$ & 2.29e-04  & 3.41 & 6.00e-03  & 1.97 \\ 
$2^6$ & $2^4$ & 1.04e-04  & 1.14 & 3.48e-03  & 0.63 \\  
$2^7$ & $2^5$ & 5.09e-05  & 1.03 & 2.48e-03  & 0.63 \\  
$2^8$ & $2^6$ & 2.54e-05  & 1.00 & 1.62e-03  & 0.61 \\ \hline
\end{tabular}
\caption{For the data set A, dependence on the grid sizes $M$ and $N$ of the errors and the rate of convergence of continuous piecewise linear approximations of the state and the parameter.}
\label{tab:fractional-kernel}
\end{center}
\end{table}
\begin{table}[t]
\begin{center}
\begin{tabular}{| l | l | c | c | c | c | }
\multicolumn{6}{c}{$\varepsilon=2^{-9}$}            \\ \hline
$N$   &	$M$   & $e_{u,2}$ & rate & $e_\vt^*$ & rate \\ \hline
$2^4$ &	$2^2$ & 1.41e-04  & -    & 1.28e-02  & -    \\
$2^5$ &	$2^3$ & 3.24e-05  & 2.12 & 3.16e-03  & 2.02 \\
$2^6$ &	$2^4$ & 7.66e-06  & 2.08 & 7.79e-04  & 2.02 \\
$2^7$ &	$2^5$ & 1.90e-06  & 2.01 & 1.94e-04  & 2.01 \\ \hline
\end{tabular}
\hspace{.2cm}
\begin{tabular}{| l | l | c | c | c | c | }
\multicolumn{6}{c}{$\varepsilon=2^{-4}$}            \\ \hline
$N$   & $M$   & $e_{u,2}$ & rate & $e_\vt^*$ & rate \\ \hline
$2^4$ & $2^2$ & 1.66e-04  & -    & 1.31e-02  & -    \\	
$2^5$ & $2^3$ & 3.55e-05  & 2.22 & 3.10e-03  & 2.08 \\			
$2^6$ & $2^4$ & 8.01e-06  & 2.14 & 7.59e-04  & 2.03 \\
$2^7$ & $2^5$ & 1.92e-06  & 2.06 & 1.90e-04  & 2.00 \\ \hline
\end{tabular}
\caption{For the data set B, dependence on the grid sizes $M$ and $N$ of the errors and the rate of convergence of continuous piecewise linear approximations of the state and the parameter. Results are provided for two choices of interaction radius $\varepsilon$.}
\label{tab:Hconv_exact}
\end{center}
\end{table} 
\begin{figure}[t]
\begin{center}
\begin{tabular}{cc}
\includegraphics[width=0.45\textwidth]{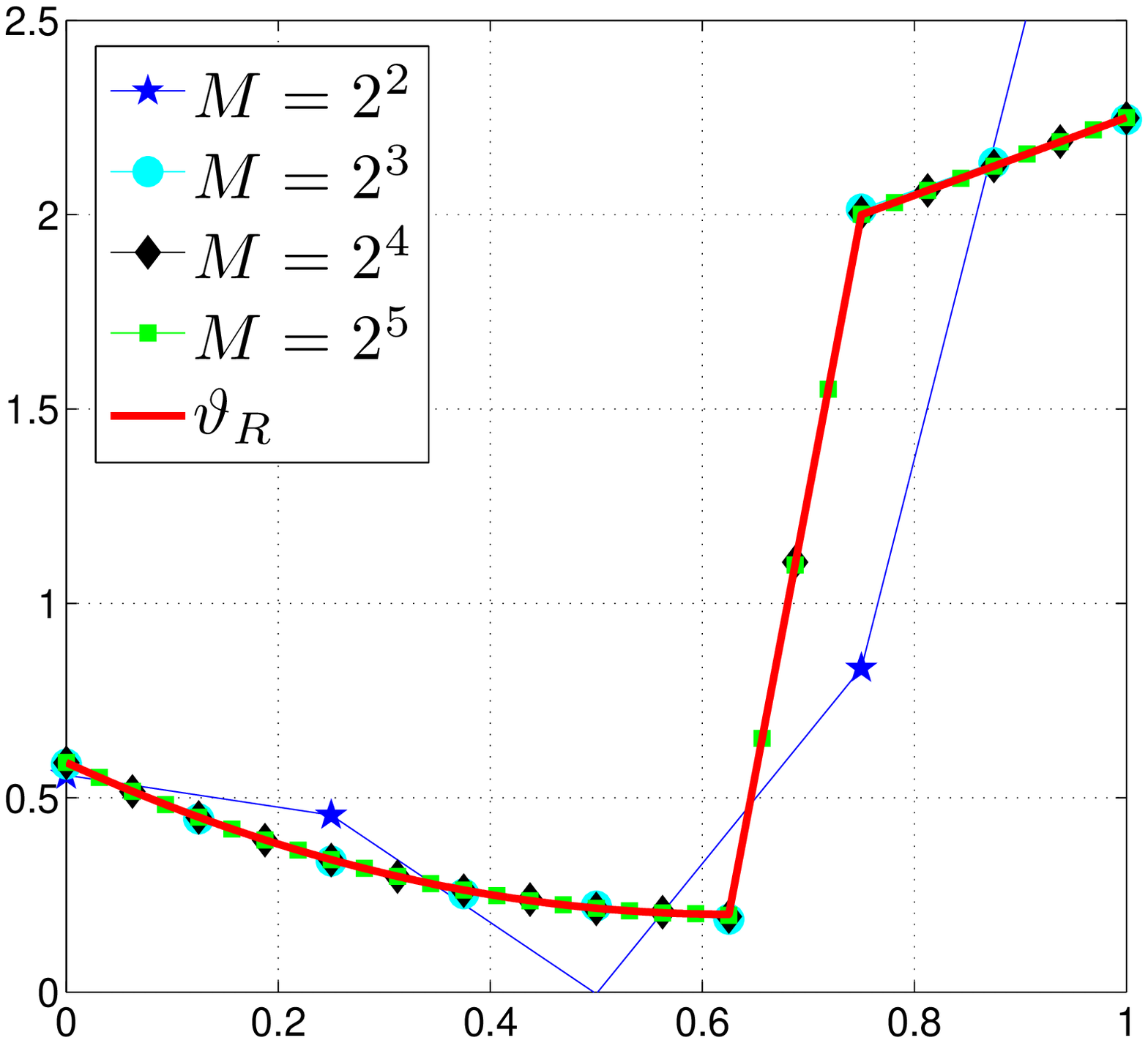} &
\includegraphics[width=0.45\textwidth]{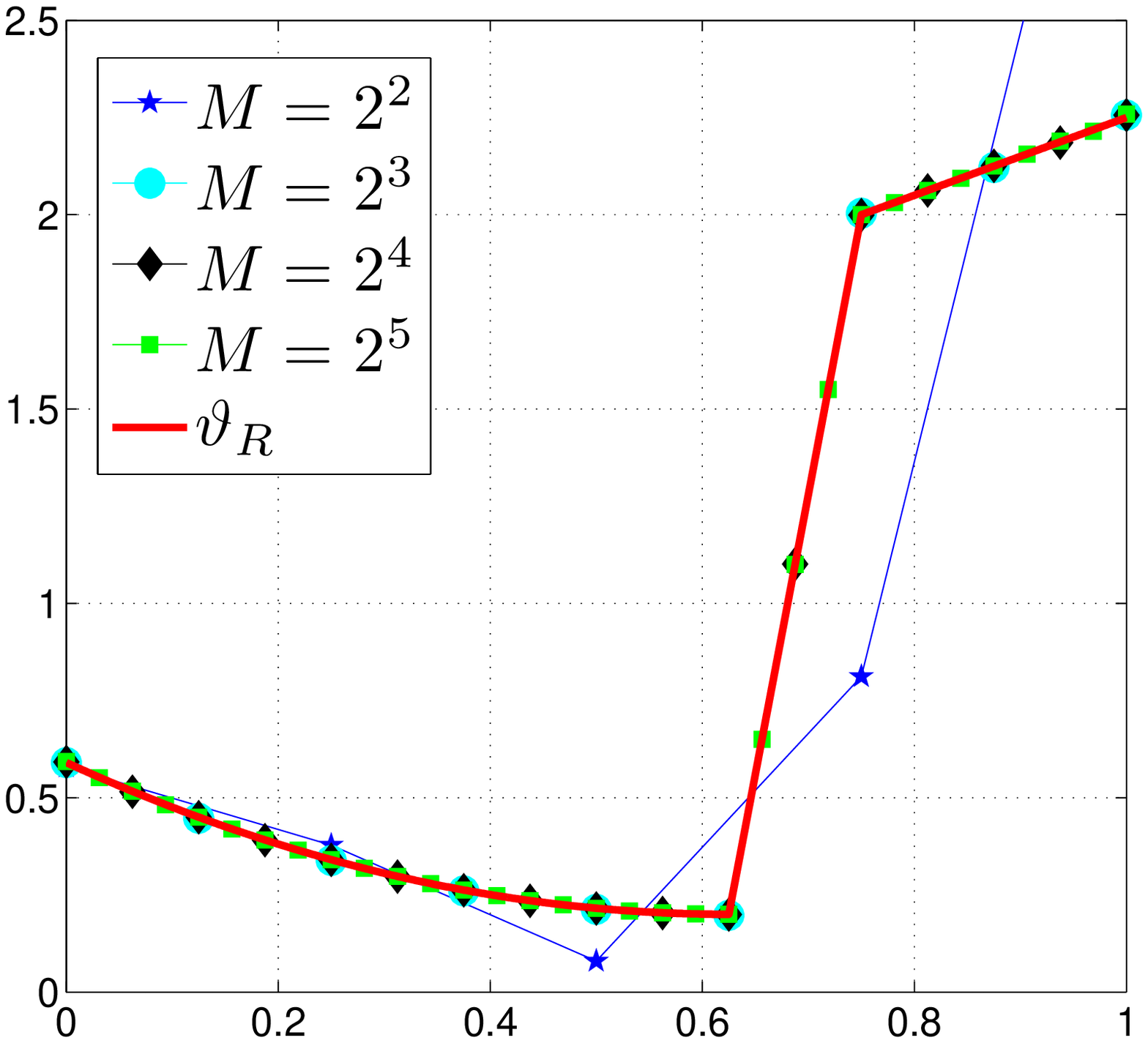} 
\end{tabular}
\end{center}
\caption{For the data set C and for different grid sizes $M$, $\vt_M^*$ and $\vt_{\rm C}$ for $\varepsilon=2^{-9}$ (left) and $2^{-4}$ (right).}
\label{fig:tDNM_DG}
\end{figure}
\begin{figure}[t]
\begin{center}
\begin{tabular}{cc}
\includegraphics[width=0.45\textwidth]{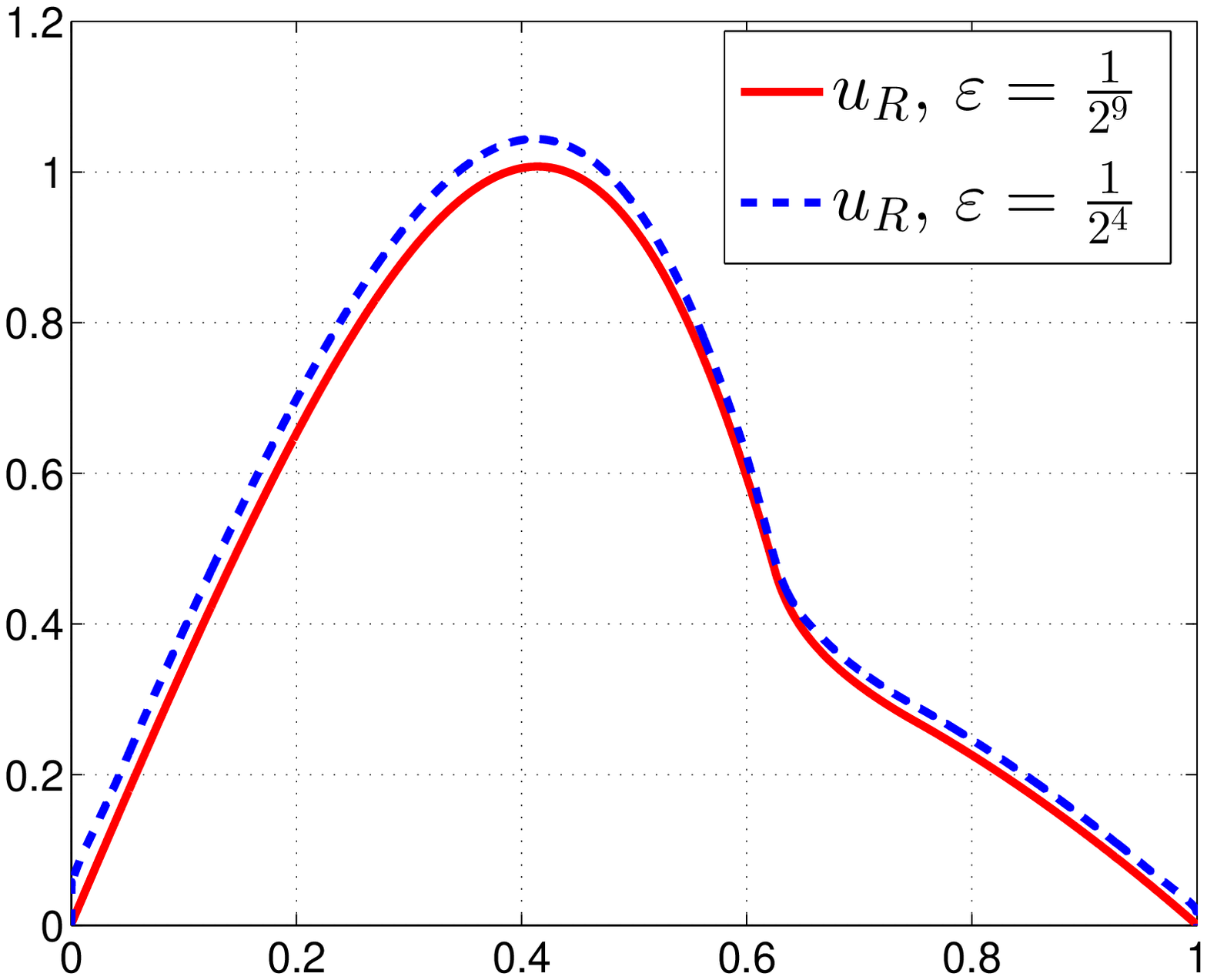} 
\includegraphics[width=0.45\textwidth]{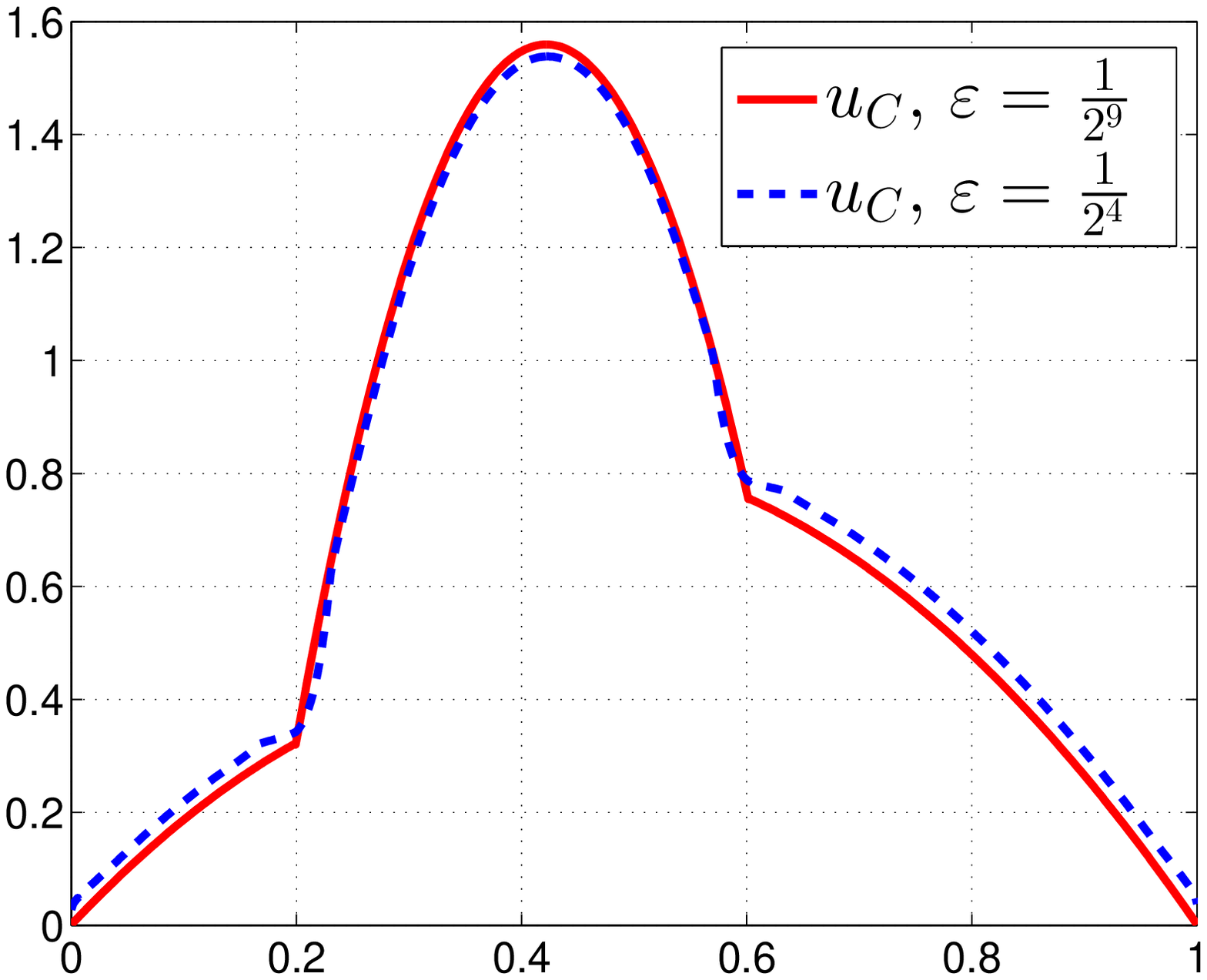} 
\end{tabular}
\end{center}
\caption{For $N=2^{12}$ and two values of $\varepsilon$, $u_{\rm C}$ on the left and $u_{\rm D}$ on the right.}
\label{fig:uREF_DG}
\end{figure}

For the data set D, $u_{\rm D}$ is determined using $N=2^{12}$ and $\beta=5 \cdot 10^{-4}$ in \eqref{regularization}. Here, $\vt_{\rm D}$ is discountinuous and for the same reasons as in case C the errors and the convergence rates are not reported. In Figure \ref{fig:tconst} we report $\vtn$ and $\vt_{\rm D}$ for several values of $M$ and we observe the convergence of the approximate optimal parameters to $\vt_{\rm D}$ as the grid is refined. In Figure \ref{fig:uREF_DG} (right) we report $u_{\rm D}$. 
\begin{figure}[t]
\begin{center}
\begin{tabular}{cc}
\includegraphics[width=0.45\textwidth]{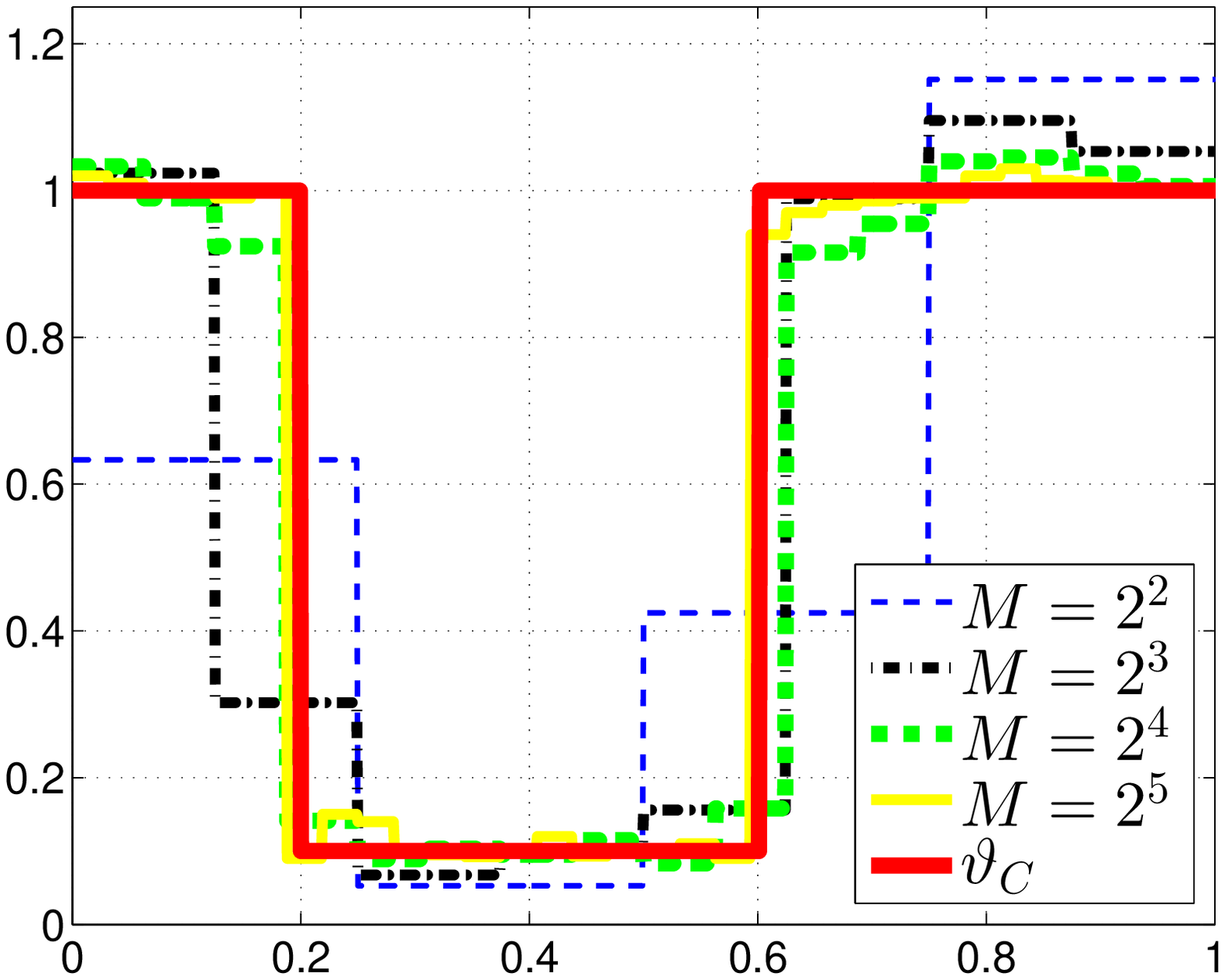} &
\includegraphics[width=0.45\textwidth]{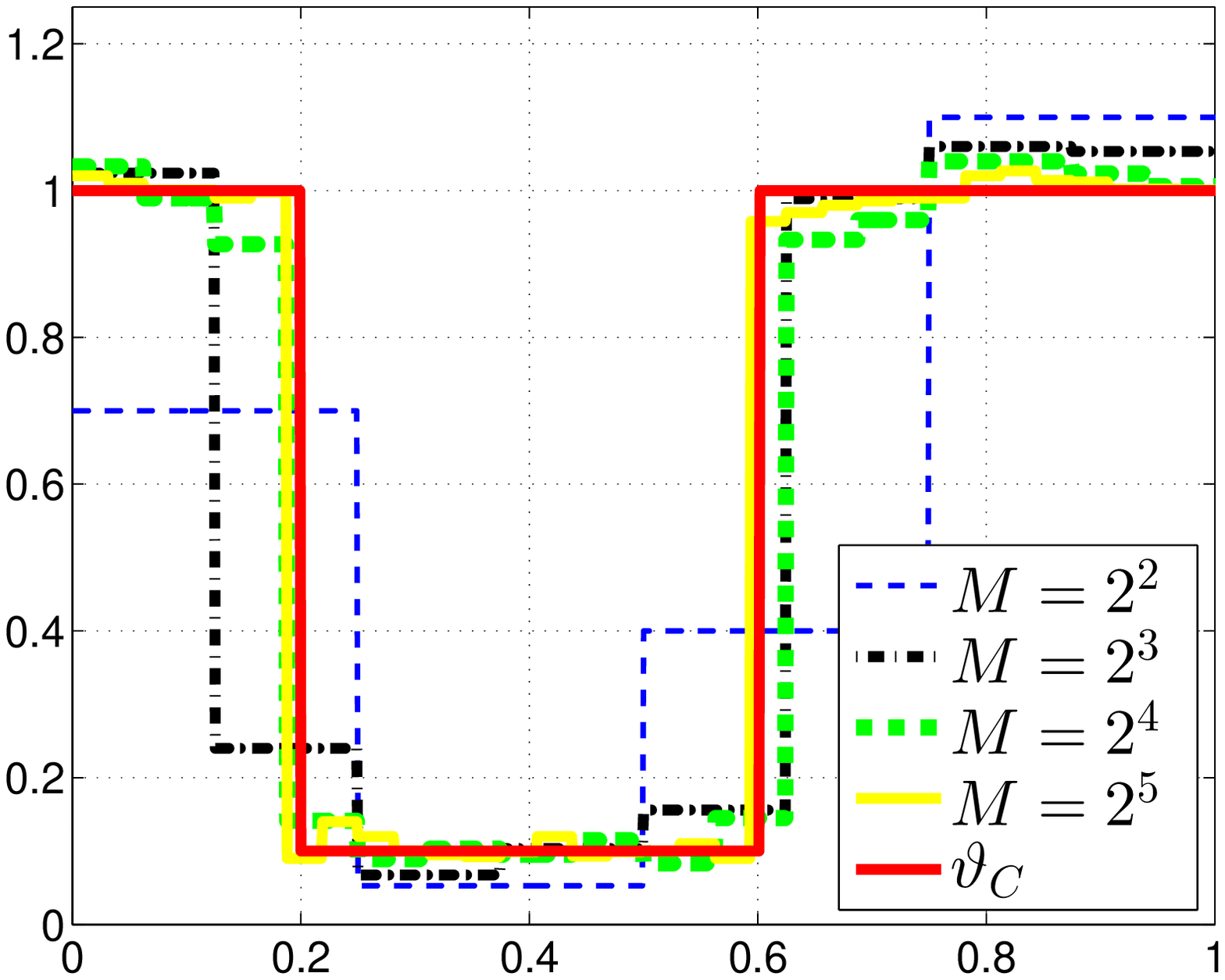} 
\end{tabular}
\end{center}
\caption{For the data set D and for different grid sizes $M$, $\vt_M^*$ and $\vt_{\rm C}$ for $\varepsilon=2^{-9}$ (left) and $2^{-4}$ (right).}
\label{fig:tconst}
\end{figure}
}
\section{Conclusions}\label{conclusion} 
The estimation of input parameters of mathematical models is an important issue also in the field of nonlocal equations. In this paper we show that the identification problem for nonlocal diffusion  equations can be treated in a way similar to the classical control theory for PDEs. Furthermore, we show that with nonlocal models we are able to estimate in an accurate way non-smooth and discontinuous input parameters. Using the nonlocal vector calculus we provide an analytical and numerical framework for the identification problem; we formulate the estimation of the diffusion parameter as an optimal control problem for which we prove the existence of an optimal solution. Using approaches similar to those for classical, local, problems we provide a Galerkin finite dimensional formulation and derive a priori error estimates for the approximate state and control variables; these theoretical results are illustrated in one-dimensional numerical examples which provide the basis for realistic simulations. Our preliminary theoretical and numerical analysis is conducted on simplified problems; obvious extensions of this work include the generalization from the scalar to the tensor diffusivity and the simulation of two- and three-dimensional problems. The former consists in estimating the coefficient $\bthe$, a second order symmetric positive definite tensor representing the diffusivity. In this case, the weak formulation of the state equation is 
\begin{displaymath}
\hspace{-2cm}
\begin{array}{ll}
   &\mbox{\em given $f(\xb)$, $g(\xb)$ and $\bthe(\xb,\yb)$, seek $u(\xb)$ such that}
   \\&\mbox{\em $u=g$ for $\xb\in\omgc$ and}
   \\&  \displaystyle \int_\omgsc\int_\omgsc \mcG u\cdot(\bthe\mcG v)\,d\yb d\xb =
   \int_\omg fv\,d\xb
   \qquad\forall\,v\in V_c(\omgsc).
\end{array}
\end{displaymath}
The latter involves more complex numerical schemes and time consuming numerical experiments \cite{chgu:11}.

An important follow-up to our work is to consider parameters affected by uncertainty. Our approach to the estimation problem is only deterministic, this is a limitation. In a stochastic approach to the solution of nonlocal problems we assume the input parameters to be random fields described by a probability distribution, i.e., $\vt = \vt(\xb,\yb, \omega)$ where $\xb$ and $\yb$ are (deterministic) points inside of $\omg$ and the random variable $\omega$ indicates that the value of the diffusion parameter at any point is drawn randomly. The same definition holds for $g(\xb,\omega)$ and $f(\xb,\omega)$; the solution of the nonlocal problem $u(\xb,\omega)$ is a random field itself. The nonlocal stochastic problem is then formulated as  
\begin{equation}\label{weakstoch}
\hspace{-2cm}\begin{array}{ll}
   &\mbox{\em given $f(\xb,\omega)$, $g(\xb,\omega)$, and $\vt(\xb,\yb,\omega)$, seek $u(\xb,\omega)$ such that}
   \\&\mbox{\em $u=g$ for $\xb\in\omgc$ and}\; \mcD(\vt\, \mcG u) = f.
\end{array}
\end{equation}
However, in practice one usually does not know much about the statistics of the input variables; the only known quantities might be the maximum and minimum values or the mean and covariance. Thus, as for the local case, a fundamental issue in nonlocal stochastic equations is the estimation of the distributions of the input parameters. In the classical framework this problem is known as {\it model calibration}. We formulate the stochastic parameter estimation problem for a nonlocal diffusion equation as the standard estimation problem for PDEs.
\begin{displaymath}
\hspace{-2cm}\begin{minipage}{4.9in}
{\em Given the random fields $f(\xb,\omega)$, $g(\xb,\omega)$, and the target function $\widehat u(\xb)$, seek an optimal control $\vt^\ast(\xb,\yb,\omega)$ and an optimal state $u^\ast(\xb,\omega)$ such that $J(u,\vt)$ is minimized, subject to $u$ and $\vt$ satisfying \eqref{weakstoch},} 
\end{minipage}
\end{displaymath}
where
\begin{displaymath}
J(u,\vt) = \mathbb{E}[\|u(\xb,\omega)-\widehat{u}(\xb)\|^2_{L^2(\omg)}].
\end{displaymath}
Here, $\mathbb{E}$ represents the expected value. We might need to add a regularization term to prevent potential ill-posedness or ill-conditioning of the analytical or numerical problem.\\


\end{document}